\documentclass[12pt]{article}
\usepackage{amsmath}
\usepackage{amsthm}
\usepackage{amssymb}
\usepackage{amsfonts}
\usepackage{epsfig}

\makeatletter
\newfont{\footsc}{cmcsc10 at 8truept}
\newfont{\footbf}{cmbx10 at 8truept}
\newfont{\footrm}{cmr10 at 10truept}
\makeatother 

\title{Prefix exchanging and pattern avoidance by involutions}
\author{Aaron D.\ Jaggard\thanks{This work is drawn from the
author's Ph.D. dissertation which was written at the University of
Pennsylvania under the supervision of Herbert S.\ Wilf. The author
was partially supported by the DoD University Research Initiative
(URI) program administered by the ONR under grant
N00014-01-1-0795; the presentation of this work at the Permutation
Patterns 2003 conference was partially supported by Penn GAPSA and
the New Zealand Institute for Mathematics and its Applications.}\\
\small Department of Mathematics\\[-0.8ex]
\small Tulane University\\[-0.8ex]
\small New Orleans, LA 70118 USA\\[-0.8ex]
\small \texttt{adj@math.tulane.edu}}
\date{\small MR Subject Classifications: 05A05, 05A15}

\def\S #1{\mathcal{S}_{#1}}

\newcommand{\iav}[2]{\mathcal{I}_{#1}(#2)}

\def\seq {\sim_\mathcal{S}}
\def\ieq {\sim_\mathcal{I}}
\def\sq #1{SQ_{#1}}
\def\sqn {\sq{n}}

\newcommand{\f}[3]{f^{#1}_{{#2},{#3}}}
\newcommand{\g}[3]{g^{#1}_{{#2},{#3}}}

\DeclareMathOperator{\evac}{evac}

\theoremstyle{plain}
\newtheorem{theorem}{Theorem}[section]
\newtheorem{lemma}[theorem]{Lemma}
\newtheorem{prop}[theorem]{Proposition}
\newtheorem{cor}[theorem]{Corollary}

\theoremstyle{definition}
\newtheorem{defn}[theorem]{Definition}
\newtheorem*{q}{Question}
\newtheorem{conj}[theorem]{Conjecture}

\newtheorem{ex}[theorem]{Example}

\theoremstyle{remark}
\newtheorem{remark}[theorem]{Remark}

\theoremstyle{plain}

\begin{document}

\maketitle

\begin{abstract}
Let $I_n(\pi)$ denote the number of involutions in the symmetric
group $\S{n}$ which avoid the permutation $\pi$.  We say that two
permutations $\alpha,\beta\in\S{j}$ \textit{may be exchanged} if
for every $n$, $k$, and ordering $\tau$ of $j+1,\ldots,k$, we have
$I_n(\alpha\tau)=I_n(\beta\tau)$.  Here we prove that $12$ and
$21$ may be exchanged and that $123$ and $321$ may be exchanged.
The ability to exchange $123$ and $321$ implies a conjecture of
Guibert, thus completing the classification of $\S{4}$ with
respect to pattern avoidance by involutions; both of these results
also have consequences for longer patterns.

Pattern avoidance by involutions may be generalized to rook
placements on Ferrers boards which satisfy certain symmetry
conditions.  Here we provide sufficient conditions for the
corresponding generalization of the ability to exchange two
prefixes and show that these conditions are satisfied by $12$ and
$21$ and by $123$ and $321$. Our results and approach parallel
work by Babson and West on analogous problems for pattern
avoidance by general (not necessarily involutive) permutations,
with some modifications required by the symmetry of the current
problem.
\end{abstract}

\section{Introduction}
\label{sec:intro} The \textit{pattern} of a sequence
$w=w_1w_2\ldots w_k$ of $k$ distinct letters is the
order-preserving relabelling of $w$ with $[k]=\{1,2,\ldots,k\}$.
Given a permutation $\pi=\pi_1\pi_2\ldots\pi_n$ in the symmetric
group $\S{n}$, we say that $\pi$ \textit{avoids} the pattern
$\sigma=\sigma_1\sigma_2\ldots\sigma_k\in\S{k}$ if there is no
subsequence $\pi_{i_1}\ldots\pi_{i_k}$, $i_1<\cdots<i_k$, whose
pattern is $\sigma$.

Let $\iav{n}{\sigma}$ denote the number of involutions in $\S{n}$
(permutations whose square is the identity permutation) which
avoid the pattern $\sigma$, and write $\sigma\ieq\sigma'$ if for
every $n$, $\iav{n}{\sigma}=\iav{n}{\sigma'}$.  (In this case we
also say that $\sigma$ and $\sigma'$ are
\textit{$\ieq$-equivalent\/}.)\ \ For $\alpha,\beta\in\S{j}$, we
say that the prefixes $\alpha$ and $\beta$ \textit{may be
exchanged} if for every $k\geq j$ and ordering
$\tau=\tau_{j+1}\tau_{j+2}\ldots\tau_{k}$ of $[k]\setminus[j]$,
the patterns
$\alpha_1\ldots\alpha_j\tau_{j+1}\tau_{j+2}\ldots\tau_{k}$ and
$\beta_1\ldots\beta_j\tau_{j+1}\tau_{j+2}\ldots\tau_{k}$ are
$\ieq$-equivalent.

Our work here implies the following corollaries about the ability
to exchange certain prefixes.  These results and the techniques we
use throughout this paper closely parallel work by Babson and
West~\cite{Babson2000} on pattern avoidance by general
permutations (without the restriction to involutions).
\newtheorem*{ex12cor}{Corollary~\ref{cor:ex12}}
\begin{ex12cor}
The prefixes $12$ and $21$ may be exchanged.
\end{ex12cor}
\newtheorem*{ex123cor}{Corollary~\ref{cor:ex123}}
\begin{ex123cor}
The prefixes $123$ and $321$ may be exchanged.
\end{ex123cor}
\noindent Corollary~\ref{cor:ex123} implies an affirmative answer
to a conjecture of Guibert (that $1234\ieq 1432$) and thus
completes the classification of $\S{4}$ according to
$\ieq$-equivalence.  Corollaries~\ref{cor:ex12}
and~\ref{cor:ex123} also imply $\ieq$-equivalences for patterns of
length greater than $4$; we discuss these in some detail for
patterns in $\S{5}$.

These corollaries follow from the sufficient conditions for
exchanging prefixes given by Corollary~\ref{cor:preex}. Recent
work by Stankova and West~\cite{Stankova2002} and
Reifegerste~\cite{Reifegerste2003} on different aspects of pattern
avoidance by general permutations suggests the generalization of
Corollary~\ref{cor:preex} given by Theorem~\ref{thm:geninvst}
below.  In order to state this theorem, we need the following
definitions.
\begin{defn}
Given a (Ferrers) shape $\lambda$, a \textit{placement on
$\lambda$} is an assignment of dots to some of the boxes in
$\lambda$ such that no row or column contains more than one dot.
We call a placement on $\lambda$ \textit{full} if each row and
column of $\lambda$ contains exactly $1$ dot.  We define the
\textit{transpose} of a placement to be the placement which has a
dot in box $(i,j)$ if and only if the original placement had a dot
in box $(j,i)$.  We call a placement on a shape $\lambda$
\textit{symmetric} if the transpose of the placement is the
original placement.
\end{defn}
\noindent The transpose of a placement on a shape $\lambda$ is a
placement on the conjugate shape $\lambda'$ of $\lambda$.  We use
`self-conjugate' to describe the symmetry of shapes and
`symmetric' to describe the symmetry of placements on shapes; our
work makes use of symmetric placements on self-conjugate shapes.

\begin{figure}[ht]
\begin{center}
    \includegraphics[height=2cm]{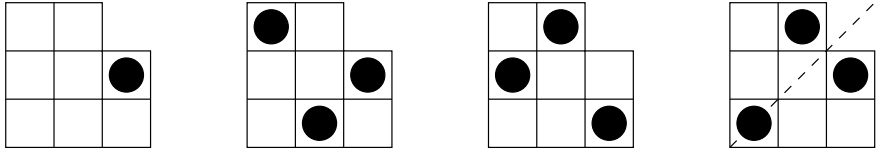}
\end{center}
\caption{Four placements on the self-conjugate shape
$(3,3,2)$.\label{fig:termex}}
\end{figure}

Figure~\ref{fig:termex} shows four placements on the
self-conjugate shape $\lambda=(3,3,2)$.  The placement on the far
left has one dot and is not full.  The placement on the center
left of the figure is full but not symmetric; its transpose is
shown at the center right of this figure.  Finally, the placement
on the far right is a symmetric full placement, with the dashed
line indicating the symmetry of the placement.

\begin{figure}[h]
\begin{center}
    \includegraphics[height=2cm]{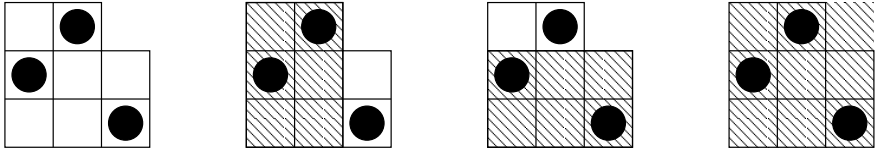}
\end{center}
\caption{A placement on $(3,3,2)$ which contains the patterns $12$
and $21$ but not the pattern $231$}\label{fig:shcon:ex}
\end{figure}

Pattern containment can be generalized to placements on shapes (as
in, \textit{e.g.\/}, \cite{Babson2000}) as follows.
\begin{defn}
A placement on a shape $\lambda$ \textit{contains the pattern
$\sigma$} if there is a set $\{(x_i,y_i)\}_{i\in[j]}$ of $j$ dots
in the placement which are in the same relationship as the values
of $\sigma$ (\textit{i.e.\/}, $x_1<\cdots<x_j$ and the pattern of
$y_1\ldots y_j$ is $\sigma$) and which are bounded by a
rectangular subshape of $\lambda$. \label{def:shapecontain}
\end{defn}
\begin{ex}
Figure~\ref{fig:shcon:ex} shows a placement on the shape $(3,3,2)$
which contains the patterns $12$ and $21$; dots whose heights form
these patterns are bounded by the shaded rectangular subshapes of
$(3,3,2)$ indicated in the center left and right of
Figure~\ref{fig:shcon:ex}.  This placement does not contain the
pattern $231$ because, although the heights of the dots in the
placement form the pattern $231$, the smallest rectangular shape
(shaded, far right) which bounds all three of these dots is not a
subshape of $(3,3,2)$.
\end{ex}

With these definitions in hand, we may state the most general
theorem that we prove here.
\newtheorem*{genthm}{Theorem~\ref{thm:geninvst}}
\begin{genthm}
Let $\lambda_{sym}(T)$ be the number of symmetric full placements
on the shape $\lambda$ which avoid all of the patterns in the set
$T$.  Let $\alpha$ and $\beta$ be involutions in $\S{j}$.  Let
$T_\alpha$ be a set of patterns, each of which begins with the
prefix $\alpha$, and $T_\beta$ be the set of patterns obtained by
replacing in each pattern in $T_\alpha$ the prefix $\alpha$ with
the prefix $\beta$.  If for every self-conjugate shape $\lambda$
$\lambda_{sym}(\{\alpha\}) = \lambda_{sym}(\{\beta\})$, then for
every self-conjugate shape $\mu$
\begin{equation*}
\mu_{sym}(T_\alpha)=\mu_{sym}(T_\beta).
\end{equation*}
\end{genthm}
\noindent  Here we also prove that the conditions on $\alpha$ and
$\beta$ in Theorem~\ref{thm:geninvst} are satisfied by the
patterns $12$ and $21$ (Theorem~\ref{thm:mu12}) and by $123$ and
$321$ (Theorem~\ref{thm:mu123}). Corollaries~\ref{cor:ex12}
and~\ref{cor:ex123} then follow.

Section~\ref{sec:back} reviews the work mentioned above and other
relevant literature and gives some additional basic definitions.
Section~\ref{sec:gen} contains some general theorems related to
involutions and patterns.  In Sections~\ref{sec:12tau}
and~\ref{sec:123tau} we show that we can apply this general
machinery to the prefixes $12$ and $21$ and then to $123$ and
$321$. Finally, in Section~\ref{sec:conc} we discuss some
$\ieq$-equivalences implied by our work as well as some
interesting open questions.

\section{Background}
\label{sec:back}

\subsection{General preliminaries}

\begin{figure}[h]
\begin{center}
    \includegraphics[width=4in]{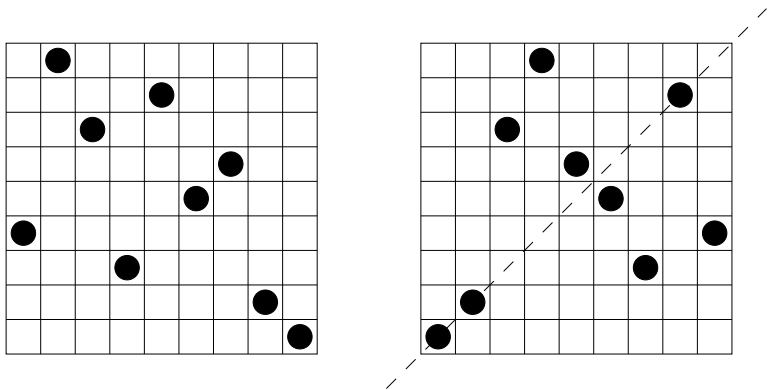}
\end{center}
  \caption{The graphs of the non-involution $497385621$ (left)
  and the involution $127965384$ (right).}\label{fig:permgraphs}
\end{figure}

We make use of the following representation of a permutation.
\begin{defn} The \textit{graph of a permutation}
$\pi\in\S{n}$ is an $n\times n$ array of boxes with dots in
exactly the set of boxes $\{(i,\pi(i))\}_{i\in[n]}$.
\label{def:graph}
\end{defn}
\noindent The graph of $\pi^{-1}$ has a dot in the box $(x,y)$ if
and only if the graph of $\pi$ has a dot in the box $(y,x)$.  We
coordinatize the graphs of permutations from the bottom left
corner, so a permutation is an involution if and only if its graph
is symmetric about the diagonal connecting its bottom left and top
right corners.  The graphs of $497385621$ (a non-involution) and
$127965384$ (an involution) are shown in
Figure~\ref{fig:permgraphs}, with the dashed line indicating the
symmetry which characterizes graphs of involutions. Denoting by
$\sqn$ the $n\times n$ square shape, the graphs of involutions in
$\S{n}$ are exactly the symmetric full placements on $\sqn$.

We will make use of the Robinson-Schensted-Knuth (RSK) algorithm,
which is treated in both Chapter~7 and Appendix~1 of~\cite{EC2}.
The RSK algorithm gives a bijection between permutations in
$\S{n}$ and pairs $(P,Q)$ of standard Young tableaux such that the
shape of $P$ is that of $Q$ and this common shape has $n$ boxes.
If $\pi\leftrightarrow (P,Q)$, then $\pi^{-1}\leftrightarrow
(Q,P)$, so this gives a bijection between $n$-involutions and
single tableaux of size $n$.

The \textit{Sch\"utzenberger involution}, or \textit{evacuation},
is an operation on tableaux.  Given a tableau $Q$, it produces a
tableau $\evac(Q)$ of the same shape as $Q$ and such that
$\evac(\evac(Q))=Q$.  A complete development of this operation is
given in Appendix~1 of~\cite{EC2}.  We note here the following
property, due to Sch\"utzenberger, which is given as
Corollary~A1.2.11 in~\cite{EC2}.
\begin{prop}[Sch\"utzenberger \cite{Schutzenberger1963}]
Let $w=w_1\ldots w_n\leftrightarrow (P,Q)$.  Then
\begin{equation*}
w_n\ldots w_1 \leftrightarrow (P^t,\evac(Q)^t)
\end{equation*}
where $P^t$ denotes the transpose of the tableau $P$.
\label{prop:evacrev}
\end{prop}

\subsection{Pattern avoidance background}

As for pattern avoidance by permutations in general, some
$\ieq$-equivalences follow from symmetry considerations. Four of
the symmetries of the square preserve the symmetry which
characterizes the graphs of involutions.  The images of a pattern
$\tau$ under these symmetries are patterns which are trivially
$\ieq$-equivalent to $\tau$; these patterns form the
\textit{(involution) symmetry class\/} of $\tau$.  For
$\tau\in\S{k}$ these patterns are $\tau$, $\tau^{-1}$, the
reversed complement
$\tau_{rc}=(k+1-\tau_k)\ldots(k+1-\tau_2)(k+1-\tau_1)$ of $\tau$,
and $(\tau_{rc})^{-1}$.  Since we cannot use all $8$ of the
symmetries of the square, each symmetry class which arises in
pattern avoidance by general permutations may split into $2$
involution symmetry classes.  We refer to the $\ieq$-equivalence
classes as \textit{(involution) cardinality classes\/}; for
pattern avoidance by general permutations, the cardinality classes
are usually referred to as \textit{Wilf classes\/}.  Unless
otherwise stated, we take `symmetry' and `cardinality' classes to
be with respect to $\ieq$, and use `$\seq$-' to indicate
equivalence with respect to pattern avoidance by general
permutations.

In their well-known paper~\cite{Simion1985}, Simion and Schmidt
found the cardinality classes of $\S{3}$ and proved the following
propositions.
\begin{prop}[Simion and Schmidt \cite{Simion1985}]
For $\tau\in\{123,132,213,321\}$ and $n\geq 1$,
\begin{equation*}
\iav{n}{\tau}=\binom{n}{[n/2]}.
\end{equation*}
\end{prop}
\begin{prop}[Simion and Schmidt \cite{Simion1985}]
For $\tau\in\{231,312\}$ and $n\geq 1$,
\begin{equation*}
\iav{n}{\tau}=2^{n-1}.
\end{equation*}
\end{prop}
\noindent Comparing this to the classic result that $\S{3}$
contains a single Wilf class, we see that passing from symmetry to
cardinality classes does not repair all of the breaks in
$\seq$-symmetry classes caused by considering pattern avoidance by
involutions instead of general permutations.

Many of the sequences $\{\iav{n}{\tau}\}$ which are known are for
$\tau=12\ldots k$, in which case the sequence counts the number of
standard tableaux of size $n$ with at most $k-1$ columns.  A
theorem of Regev covers $k=4$ as follows.
\begin{prop}[Regev \cite{Regev1981}]
\begin{equation*}
\iav{n}{1234}=\sum_{i=0}^{\lfloor n/2\rfloor} \binom{n}{2i}
\binom{2i}{i}\frac{1}{i+1},
\end{equation*}
\textit{i.e.\/}, the $n^\mathrm{th}$ Motzkin number $M_n$.
\end{prop}
\noindent Regev also gave the following expression for the
asymptotic value of $\iav{n}{12\ldots k}$.
\begin{theorem}[Regev {\cite{Regev1981}}]
\begin{equation*}
\iav{n}{12\ldots k(k+1)}\sim k^n
\left(\frac{k}{n}\right)^{k(k-1)/4} \frac{1}{k!}
\Gamma(\frac{3}{2})^{-k} \prod_{j=1}^k \Gamma(1+\frac{j}{2})
\end{equation*}
\end{theorem}
\noindent Gouyou-Beauchamps studied Young tableaux of bounded
height in~\cite{G-B1989} and obtained exact results for $k=5$ and
$6$.
\begin{prop}[Gouyou-Beauchamps \cite{G-B1989}]
\begin{equation*}
\iav{n}{12345}=\begin{cases}
C_kC_k, &n=2k-1\\
C_kC_{k+1}, &n=2k
\end{cases},
\end{equation*}
where $C_k=\frac{1}{k+1}\binom{2k}{k}$, the $k^\mathrm{th}$
Catalan number.
\end{prop}
\begin{prop}[Gouyou-Beauchamps \cite{G-B1989}]
\begin{equation*}
\iav{n}{123456} = \sum_{i=0}^{\lfloor n/2\rfloor}
\frac{3!n!(2i+2)!}{(n-2i)!i!(i+1)!(i+2)!(i+3)!}.
\end{equation*}
\end{prop}
\noindent Gessel~\cite{Gessel1990} has given a determinantal
formula for the general $\iav{n}{12\ldots k}$.

Work of Guibert and others has almost completely determined the
cardinality classes of $\S{4}$ (see~\cite{Guibert2001} for a
review of this work).  Symmetry of the RSK algorithm implies
$1234\ieq 4321$.  Guibert bijectively obtained the following
results in his thesis~\cite{GuibertThesis}.
\begin{prop}[Guibert \cite{GuibertThesis}]
\begin{equation*}
3412\ieq 4321
\end{equation*}
\end{prop}
\begin{prop}[Guibert \cite{GuibertThesis}]
\begin{equation*}
2143\ieq 1243
\end{equation*}
\end{prop}
\noindent Guibert also conjectured that both $\iav{n}{2143}$ and
$\iav{n}{1432}$ are equal to $M_n$ for $n\geq 4$ (as
$\iav{n}{1234}$ is known to be).  Guibert, Pergola, and
Pinzani~\cite{Guibert2001} affirmatively answered the first of
these conjectures.
\begin{prop}[Guibert, Pergola, Pinzani \cite{Guibert2001}]
\begin{equation*}
1234\ieq 2143
\end{equation*}
\end{prop}
\noindent In more recent work on involutions avoiding various
combinations of multiple patterns, Guibert and
Mansour~\cite{Guibert2002b} noted that the second conjecture was
still open.  We prove that conjecture as
Corollary~\ref{cor:inv:3214}.

There are various known $\seq$-equivalences between
$\seq$-symmetry classes.  Of particular interest are those which
follow from more general theorems, which we review here.  In
Sections~\ref{sec:gen}--\ref{sec:123tau} we prove the first such
general theorems for pattern avoidance by involutions.

West proved the following theorem in his thesis~\cite{WestThesis}.
\begin{theorem}[West \cite{WestThesis}]
For any $k$, any ordering $\tau=\tau_3\ldots\tau_k$ of
$[k]\setminus[2]$, and any $n$, the number of permutations in
$\S{n}$ which avoid the pattern $12\tau_3\ldots\tau_k$ equals the
number of permutations in $\S{n}$ which avoid
$21\tau_3\ldots\tau_k$. \label{thm:perm12tau}
\end{theorem}
\noindent Babson and West~\cite{Babson2000} restated the proof of
Theorem~\ref{thm:perm12tau} and then proved the following theorem.
\begin{theorem}[Babson and West \cite{Babson2000}]
For any $k$, any ordering $\tau=\tau_4\ldots\tau_k$ of
$[k]\setminus[3]$, and any $n$, the number of permutations in
$\S{n}$ which avoid the pattern $123\tau_4\ldots\tau_k$ equals the
number of permutations in $\S{n}$ which avoid
$321\tau_4\ldots\tau_k$. \label{thm:perm123tau}
\end{theorem}
Stankova and West~\cite{Stankova2002} further investigated the
property, which they called \textit{shape-Wilf-equivalence}, used
by Babson and West in their proofs of these two theorems.  Two
patterns $\alpha$ and $\beta$ are shape-Wilf-equivalent if, for
every shape $\lambda$, the number of full placements on $\lambda$
which avoid $\alpha$ equals the number which avoid $\beta$; this
implies the Wilf-equivalence of the patterns in question. Stankova
and West proved that the patterns $231\tau_4\ldots\tau_k$ and
$312\tau_4\ldots\tau_k$ are shape-Wilf-equivalent. More recently,
Backelin, West, and Xin~\cite{Backelin2003} have proved that the
patterns $12\ldots j\tau_{j+1} \ldots \tau_k$ and $j\ldots
21\tau_{j+1}\ldots\tau_k$ are shape-Wilf-equivalent.  Here we
define and use a symmetrized version of shape-Wilf-equivalence.

Finally, a recent paper by Reifegerste~\cite{Reifegerste2003}
generalizes a bijection given by Simion and Schmidt.  One
application (Corollary~9 of~\cite{Reifegerste2003}) is that a
certain set of patterns with prefix $12$ is as restrictive (with
respect to pattern avoidance by general permutations) as the set
of patterns obtained by replacing these occurrences of $12$ with
$21$; this suggests part of our most general result below.

\section{Some General Machinery}
\label{sec:gen}

In order to prove Corollaries~\ref{cor:ex12} and~\ref{cor:ex123}
we need only Corollary~\ref{cor:preex} below and some additional
lemmas. The recent work, discussed in Section~\ref{sec:back}, by
Reifegerste and by Stankova and West suggests the generalization
of Corollary~\ref{cor:preex} given by Theorem~\ref{thm:geninvst}.
\begin{theorem}
Let $\lambda_{sym}(T)$ be the number of symmetric full placements
on the shape $\lambda$ which avoid all of the patterns in the set
$T$.  Let $\alpha$ and $\beta$ be involutions in $\S{j}$.  Let
$T_\alpha$ be a set of patterns, each of which begins with the
prefix $\alpha$, and $T_\beta$ be the set of patterns obtained by
replacing in each pattern in $T_\alpha$ the prefix $\alpha$ with
the prefix $\beta$.  If for every self-conjugate shape $\lambda$
$\lambda_{sym}(\{\alpha\}) = \lambda_{sym}(\{\beta\})$, then for
every self-conjugate shape $\mu$
\begin{equation*}
\mu_{sym}(T_\alpha)=\mu_{sym}(T_\beta).
\end{equation*}
\label{thm:geninvst}
\end{theorem}

The proof of Theorem~\ref{thm:geninvst} makes use of the following
definition.
\begin{defn}
Fix positive integers $j$ and $l$, and for every $i\in[l]$ let
$\tau_i$ be an ordering of $[k_i]\setminus[j]$ for some $k_i\geq
j$.  Let $T$ be the set $\{\tau_i\}_{i\in[l]}$ and $\mu$ be a
self-conjugate shape with a symmetric full placement $P$.  We
construct the \textit{self-conjugate $T$-shape of $(\mu,P)$},
denoted $\lambda_T(\mu,P)$, as follows; Example~\ref{ex:ltmp} and
Figure~\ref{fig:inv:ltmp} below illustrate this procedure.

Take all boxes $(x,y)$ and $(y,x)$ in $\mu$ such that $(x,y)$ is
strictly southwest of an occurrence of the pattern of some
$\tau_i\in T$ (\textit{i.e.\/}, for which there is a set of $k_i -
j$ dots, contained within a rectangular subshape of $\mu$, whose
heights have pattern $\tau_i$ and which are all above and to the
right of $(x,y)$.)\ \ This set of boxes forms a self-conjugate
shape, since it contains $(x,y)$ iff it contains $(y,x)$, on which
there is a (not necessarily full) symmetric placement obtained by
restricting $P$ to this shape. Delete the rows and columns of this
shape which do not contain a dot to obtain the self-conjugate
shape $\lambda_T(\mu,P)$.  The deletion of empty rows and columns
yields a symmetric full placement on $\lambda_T(\mu,P)$; we call
this the \textit{placement on $\lambda_T(\mu,P)$ induced by
$P$\/}.\label{def:pre:ltmp}
\end{defn}

\begin{ex}
We view the graph of $127965384$, shown in the right part of
Figure~\ref{fig:permgraphs}, as a placement $P$ on $\mu=\sq{9}$
and let $T=\{54\}$.  The left of Figure~\ref{fig:inv:ltmp} shows
this graph with shading added to those boxes which are southwest
of some pair of dots whose pattern is $21$ (the pattern of $54\in
T$) or which are the reflection of such a box across the diagonal
of symmetry. Removing the empty rows and columns from this shaded
shape, we obtain $\lambda_{\{54\}}(\sq{9},P)$ and the placement on
$\lambda_{\{54\}}(\sq{9},P)$ induced by $P$; these are shown at
the far right.  The center right of Figure~\ref{fig:inv:ltmp}
shows the graph of $127965384$ with the boxes corresponding to
$\lambda_{\{54\}}(\sq{9},P)$ crossed out. \label{ex:ltmp}
\end{ex}
\begin{figure}[ht]
\begin{center}
  \includegraphics[width=5in]{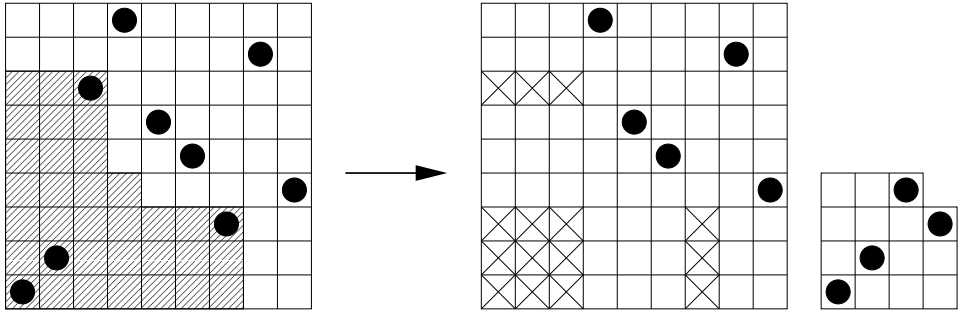}\\
\end{center}
\caption{Constructing $\lambda_T(\mu,P)$.}\label{fig:inv:ltmp}
\end{figure}

\begin{remark}
The shape of $\lambda_T(\mu,P)$ does not depend on the placement
induced on it by $P$.  For any placement $P'$ on $\mu$ which
agrees with $P$ outside of the boxes corresponding to
$\lambda_T(\mu,P)$, we have $\lambda_T(\mu,P') =
\lambda_T(\mu,P)$. \label{rem:pre:shltmp}
\end{remark}

The motivation for the definition of $\lambda_T(\mu,P)$ is that it
satisfies the following lemma.
\begin{lemma}
Fix positive integers $j$ and $l$, and for every $i\in[l]$ let
$\tau_i$ be an ordering of $[k_i]\setminus[j]$ for some $k_i\geq
j$.  Let $T$ be the set $\{\tau_i\}_{i\in[l]}$. If $P$ is a
symmetric full placement on a self-conjugate shape $\mu$ and
$\sigma$ is a $j$-involution, then $P$ contains at least one of
the patterns $\{\sigma\tau_i\}_{i\in[l]}$ if and only if the
placement on $\lambda_T(\mu,P)$ induced by $P$ contains $\sigma$.
\label{lem:inv:stscontain}
\end{lemma}
\begin{proof}
Assume $P$ contains an occurrence of $\sigma\tau_i$ in the boxes
$(x_1,y_1), \ldots, (x_{k_i},y_{k_i})$. The box $(x_j,\max_{1\leq
m\leq j}\{y_m\})$ is southwest of all of the dots in the boxes
$(x_{j+1},y_{j+1})$, $\ldots$, $(x_{j+k_i},y_{j+k_i})$, whose
pattern is that of $\tau_i$. The boxes
$(x_1,y_1),\ldots,(x_j,y_j)$, whose pattern is $\sigma$, are all
(weakly) southwest of this box and are thus contained in a
rectangular subshape of $\lambda_T(\mu,P)$.  The placement on
$\lambda_T(\mu,P)$ induced by $P$ thus contains the pattern
$\sigma$.

If the placement on $\lambda_T(\mu,P)$ induced by $P$ contains
$\sigma$, we consider the top right corner $(x,y)$ of a
rectangular subshape of $\lambda_T(\mu,P)$ that bounds a set of
dots which form an occurrence of $\sigma$.  This corresponds
(after replacing the rows and columns deleted in the construction
of $\lambda_T(\mu,P)$) to a box $(x',y')$ in $\mu$ such that
either $(x',y')$ or $(y',x')$ is strictly southwest of some set of
$k_i-j$ dots whose pattern is that of some $\tau_i\in T$ and which
are all (weakly) southwest of some box in the shape $\mu$.  If the
box $(x',y')$ satisfies this condition, then the original $j$ dots
whose pattern is $\sigma$ together with the $k_i-j$ dots just
found give a $\sigma\tau_i\in T$ pattern contained in the
placement $P$.  If it does not satisfy this condition, then by
construction of $\lambda_T(\mu,P)$ the box $(y',x')$ must do so.
The reflection of the set of $j$ dots which are southwest of
$(x',y')$ and whose pattern is $\sigma$ is a set of $j$ dots which
are southwest of $(y',x')$ and whose pattern is
$\sigma^{-1}=\sigma$. These dots combine with the $k_i-j$ dots
strictly northeast of $(y',x')$ whose pattern is $\tau_i$ (and
which are southwest of some box in $\mu$) to form the pattern
$\sigma\tau_i\in T$ contained in the placement $P$.
\end{proof}

\begin{ex}
Applying Lemma~\ref{lem:inv:stscontain} to the involution
$127965384$ from Example~\ref{ex:ltmp}, we see that $127965384$
contains $12354$ (respectively $32154$) iff the placement on
$(4^3,3)$ shown at the right of Figure~\ref{fig:inv:ltmp} contains
$123$ (respectively $321$).
\end{ex}

We now prove Theorem~\ref{thm:geninvst}.
\begin{proof}[Proof of Theorem~\ref{thm:geninvst}]
Let $T$ be obtained from $T_\alpha$ by removing the prefix
$\alpha$ from every pattern in $T_\alpha$.  (Removing $\beta$ from
every pattern in $T_\beta$ also yields $T$.)\ \ For a symmetric
full placement $P$ on $\mu$, find $\lambda_T(\mu,P)$ and note
which boxes in $\mu$ correspond to the boxes of
$\lambda_T(\mu,P)$.  Let $[P]$ be the set of symmetric full
placements on $\mu$ which agree with $P$ outside of the boxes
corresponding to $\lambda_T(\mu,P)$. By
Remark~\ref{rem:pre:shltmp}, $\lambda_T(\mu,P') =
\lambda_T(\mu,P)$ for every placement $P'\in[P]$ and the number of
placements in $[P]$ equals the number of symmetric full placements
on $\lambda_T(\mu,P)$.  By Lemma~\ref{lem:inv:stscontain}, the
number of symmetric full placements in $[P]$ which avoid
$T_\alpha$ (respectively $T_\beta$) equals the number
$(\lambda_T(\mu,P))_{sym}(\{\alpha\})$ (respectively
$(\lambda_T(\mu,P))_{sym}(\{\beta\})$) of symmetric full
placements on $\lambda_T(\mu,P)$ which avoid $\alpha$
(respectively $\beta$). By hypothesis
$(\lambda_T(\mu,P))_{sym}(\{\alpha\}) =
(\lambda_T(\mu,P))_{sym}(\{\beta\})$, and the theorem follows by
summing over all classes $[P]$.
\end{proof}

A special case of Theorem~\ref{thm:geninvst} gives sufficient
conditions for the exchange of two prefixes $\alpha$ and $\beta$.
\begin{cor} Let $\lambda_{sym}(\{\sigma\})$ be the number of
symmetric full placements on the shape $\lambda$ which avoid the
pattern $\sigma$.  Let $\alpha$ and $\beta$ be involutions in
$\S{j}$. If, for every self-conjugate shape $\lambda$ we have
$\lambda_{sym}(\{\alpha\})=\lambda_{sym}(\{\beta\})$, then the
prefixes $\alpha$ and $\beta$ may be exchanged. \label{cor:preex}
\end{cor}
\begin{proof}
For any $k$ and ordering $\tau$ of $[k]\setminus [j]$, take
$T_\alpha=\{\alpha\tau\}$, $T_\beta=\{\beta\tau\}$, and
$\mu=\sqn$.  Theorem~\ref{thm:geninvst} then implies that
$\mu_{sym}(\{\alpha\tau\}) = \mu_{sym}(\{\beta\tau\})$.  As the
symmetric full placements on $\mu$ are exactly the graphs of
$n$-involutions, we have $\iav{n}{\alpha\tau} =
\iav{n}{\beta\tau}$.  Since this does not depend on our choices of
$n$ or $\tau$, the prefixes $\alpha$ and $\beta$ may be exchanged.
\end{proof}

\section{The patterns $12$ and $21$} \label{sec:12tau}

We now show that the conditions on $\{\alpha,\beta\}$ in
Theorem~\ref{thm:geninvst} are satisfied by the patterns $12$.
\begin{lemma}
For any self-conjugate shape $\lambda$, the number
$\lambda_{sym}(\{12\})$ of symmetric full placements on $\lambda$
which avoid $12$ equals the number $\lambda_{sym}(\{21\})$ which
avoid $21$. \label{lem:comb:12fill21}
\end{lemma}
\begin{proof}
Babson and West~\cite{Babson2000} showed that if $\lambda$ has any
full placements, there is a unique full placement on $\lambda$
which avoids $12$ and a unique full placement on $\lambda$ which
avoids $21$.  If $\lambda$ is self-conjugate, the reflection of
any placement on $\lambda$ across the diagonal of symmetry gives
another placement on $\lambda$.  This placement avoids $12$ ($21$,
respectively) iff the original placement did. By the uniqueness of
the full placements which avoid $12$ and $21$, the reflected
placement must coincide with the original one and is thus
symmetric.
\end{proof}

We may thus apply Theorem~\ref{thm:geninvst} to $12$ and $21$ in
order to obtain the following result.
\begin{theorem}
Let $T_{12}$ be a set of patterns, each of which begins with the
prefix $12$, and $T_{21}$ be the set of patterns obtained by
replacing in each pattern in $T_{12}$ the prefix $12$ with the
prefix $21$.  Let $\mu_{sym}(T)$ be the number of symmetric full
placements on the shape $\mu$ which avoid every pattern in the set
$T$. For every self-conjugate shape $\mu$,
\begin{equation*}
\mu_{sym}(T_{12})=\mu_{sym}(T_{21}).
\end{equation*}
\label{thm:mu12}
\end{theorem}
As a corollary (also seen by applying Corollary~\ref{cor:preex} to
Lemma~\ref{lem:comb:12fill21}), we may exchange the prefixes $12$
and $21$.
\begin{cor}
The prefixes $12$ and $21$ may be exchanged.\hfil\qed
\label{cor:ex12}
\end{cor}
\noindent We apply Corollary~\ref{cor:ex12} to specific patterns
in Section~\ref{sec:conc}.

\section{The patterns $123$ and $321$}
\label{sec:123tau}

We now turn to the prefixes $123$ and $321$ and show that these
satisfy the conditions given in Theorem~\ref{thm:geninvst}. Our
approach closely parallels that used by Babson and West in their
work on pattern avoiding permutations that we discussed in
Section~\ref{sec:back}. We symmetrize many of their results here;
as we do so, we encounter problems with symmetric full placements
on square shapes (\textit{i.e.\/}, the graphs of involutions).  We
are able to work around these problems using various symmetry
properties of the RSK algorithm.

\begin{figure}[h]
\begin{center}
    \includegraphics[height=2cm]{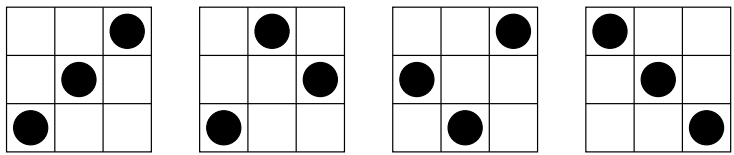}
\end{center}
\caption{The conclusion of Lemma~\ref{lem:inv:123nonsq} does not
hold for the square shape $\sq{3}$.}\label{fig:sym333}
\end{figure}

\begin{figure}[h]
\begin{center}
    \includegraphics[height=2cm]{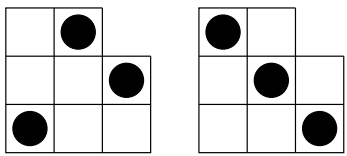}
\end{center}
\caption{Illustration of Lemma~\ref{lem:inv:123nonsq} for the
non-square shape $(3,3,2)$.}\label{fig:sym332}
\end{figure}

We start with the statement of the main lemma we use to apply
Theorem~\ref{thm:geninvst} to the patterns $123$ and $321$.  It
relates the number of $123$- and $321$-avoiding symmetric full
placements on self-conjugate shapes according to the position of
the dot in the top row.
\begin{lemma}
If $\lambda=(\lambda_1,\ldots,\lambda_k)$ is a non-square
self-conjugate shape then, for $1\leq i\leq\lambda_k$, the number
of symmetric full placements on $\lambda$ which avoid $123$ and
have a dot in $(i,\lambda_1)$ equals the number of symmetric full
placements on $\lambda$ which avoid $321$ and have a dot in
$(\lambda_k+1-i,\lambda_1)$. \label{lem:inv:123nonsq}
\end{lemma}
\noindent We defer the proof of this lemma until the end of this
section and start with an example showing that the non-square
condition in Lemma~\ref{lem:inv:123nonsq} is required.
\begin{ex}
The conclusion of Lemma~\ref{lem:inv:123nonsq} need not hold for
square shapes.  Figure~\ref{fig:sym333} shows the four symmetric
full placements on $\sq{3}$.  The three rightmost placements avoid
$123$ and have dots in $(i,3)=(2,3)$, $(3,3)$, and $(1,3)$.  The
three leftmost placements avoid $321$ and have dots in
$(4-i,3)=(3,3)$, $(2,3)$, and $(3,3)$ ($i=1$, $2$, $1$).

The two symmetric full placements on the non-square shape
$(3,3,2)$ are shown in Figure~\ref{fig:sym332}; each avoids both
$123$ and $321$.  The dots appearing at $(i,3)$ are $(2,3)$ and
$(1,3)$, while the dots at $(3-i,3)$ are $(2,3)$ and $(1,3)$
($i=1,2$).
\end{ex}

We now apply Theorem~\ref{thm:geninvst} to the patterns $123$ and
$321$.
\begin{theorem}
Let $T_{123}$ be a set of patterns, each of which begins with the
prefix $123$, and $T_{321}$ be the set of patterns obtained by
replacing in each pattern in $T_{123}$ the prefix $123$ with the
prefix $321$.  Let $\mu_{sym}(T)$ be the number of symmetric full
placements on the shape $\mu$ which avoid every pattern in the set
$T$. For every self-conjugate shape $\mu$,
\begin{equation*}
\mu_{sym}(T_{123})=\mu_{sym}(T_{321}).
\end{equation*}
\label{thm:mu123}
\end{theorem}
\begin{proof}
Summing Lemma~\ref{lem:inv:123nonsq} over $1\leq i\leq\lambda_k$,
the number of $123$-avoiding symmetric full placements on a
non-square self-conjugate shape equals the number of
$321$-avoiding such placements.  Symmetry of the RSK algorithm
gives $123\ieq 321$, so the number of symmetric full placements on
$\sqn$ which avoid $123$ equals the number which avoid $321$. We
may then apply Theorem~\ref{thm:geninvst}.
\end{proof}

As a corollary we may exchange the prefixes $123$ and $321$.
\begin{cor}
The prefixes $123$ and $321$ may be exchanged.\hfil\qed
\label{cor:ex123}
\end{cor}

\subsection*{Proving Lemma~\ref{lem:inv:123nonsq}}

The rest of this section is devoted to proving
Lemma~\ref{lem:inv:123nonsq}.  In doing so, we symmetrize the
induction used by Babson and West~\cite{Babson2000} in their proof
of an analogous lemma for pattern avoidance by general
permutations.  We start with the following lemma, a consequence of
the symmetry of the RSK algorithm, which provides additional base
cases that are needed for our symmetrized induction.  Our proof of
this lemma uses the language of $n$-involutions instead of the
(equivalent) language of symmetric full placements on $\sqn$.
\begin{lemma}
The number of symmetric full placements on $\sqn$ which avoid the
pattern $123$ and whose leftmost $i$ columns avoid $12$ equals the
number of symmetric full placements on $\sqn$ which avoid $321$
and whose rightmost $i$ columns avoid $21$. \label{lem:inv:123sq}
\end{lemma}
\begin{proof}
$123$-avoiding involutions correspond to standard Young tableaux
with at most $2$ columns.  Those which avoid $12$ in their first
$i$ entries are those whose first $i$ entries form a decreasing
subsequence; these correspond to tableaux with at most $2$ columns
and whose first column contains $1,\ldots,i$.

$321$-avoiding involutions which avoid $21$ in their last $i$
entries may be reversed to obtain $123$-avoiding permutations
which avoid $12$ in their first $i$ entries.  These correspond to
pairs $(P,Q)$ of tableaux whose common shape has at most $2$
columns and in which $Q$ contains $1,\ldots,i$ in its first
column. Proposition~\ref{prop:evacrev} shows that the pairs of
this type which correspond to the reversal of an involution are
exactly those in which $P=\evac(Q^t)^t$.
\end{proof}

We use the following lemma, which symmetrizes Lemma~2.2
of~\cite{Babson2000}, to prove Lemma~\ref{lem:inv:123nonsq} and
return to its proof to finish this section.
\begin{lemma}
Let $\lambda$ be a self-conjugate shape of length $k$, with
$i<\lambda_k$, and $1\leq j\leq\lambda_k-i$. The number of
symmetric full placements on $\lambda$ which avoid $321$ and which
avoid $21$ in the $j$ columns $i+1,\ldots,i+j$ equals the number
which avoid $321$ everywhere and $21$ in the $j$ columns
$i,\ldots,i+j-1$. \label{lem:inv:321slide}
\end{lemma}

\begin{figure}[h]
\begin{center}
    \includegraphics[width=10cm]{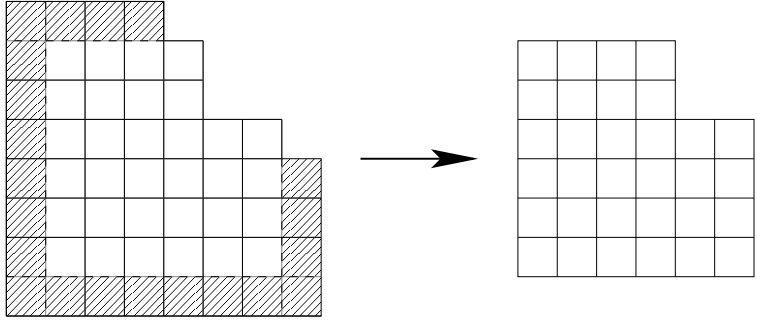}
\end{center}
\caption{A non-square shape $\lambda=(8^4,7,5^2,4)$ and the shape
$\hat{\lambda}=(6^4,4^2)$.}\label{fig:hatlambda}
\end{figure}

Our proof of Lemma~\ref{lem:inv:123nonsq} requires the following
definition.
\begin{defn}
For a non-square self-conjugate shape $\lambda$, define
$\hat{\lambda}$ to be the self-conjugate shape obtained by
deleting the leftmost and rightmost columns and the top and bottom
rows of $\lambda$. Figure~\ref{fig:hatlambda} shows the shape
$\lambda=(8^4,7,5^2,4)$ and the shape $\hat{\lambda}=(6^4,4^2)$
obtained by applying this operation.  (Note that if $\lambda$ is
non-square and non-empty and $\hat{\lambda}$ is empty, then
$\lambda = (1)$ or $(2,1)$.)
\end{defn}

We borrow notation from Babson and West to use in our symmetrized
version of their proof.
\begin{proof}[Proof of Lemma~\ref{lem:inv:123nonsq}]
For a self-conjugate shape $\lambda$ of length $k$ and $1\leq
i\leq \lambda_k$, let $\f{123}{\lambda}{i}$ be the number of
symmetric full placements on $\lambda$ which avoid $123$ and which
have a dot in $(i,k)$ (and, by symmetry, in $(k,i)$). Let
$\f{321}{\lambda}{i}$ be the number of symmetric full placements
on $\lambda$ which avoid $321$ and which have a dot in
$(\lambda_k+1-i,k)$ (and $(k,\lambda_k+1-i)$). Let
$\g{123}{\lambda}{i}$ be the number of $123$-avoiding symmetric
full placements on $\lambda$ which also avoid $12$ in the $i$
columns $1,\ldots,i$ and let $\g{321}{\lambda}{i}$ be the number
of $321$-avoiding symmetric full placements on $\lambda$ which
also avoid $21$ in the $i$ columns $\lambda_k+1-i, \ldots,
\lambda_k$. Note that the placements on $\lambda$ which avoid $12$
($21$) in certain columns also avoid $12$ ($21$) in the
same-labelled rows.  Finally, define $\g{123}{\lambda}{0}$
(respectively $\g{321}{\lambda}{0}$) to be the number of symmetric
full $123$-avoiding ($321$-avoiding) placements on $\lambda$; we
have $\g{123}{\lambda}{0}=\g{123}{\lambda}{1}$ and
$\g{321}{\lambda}{0}=\g{321}{\lambda}{1}$.

Consider a symmetric full placement on a non-square shape
$\lambda$ which is counted by $\f{123}{\lambda}{i}$.  If
$\hat{\lambda}$ is empty, $\f{123}{\lambda}{1} =
\f{321}{\lambda}{1} = 1$.  If $\hat{\lambda}$ is non-empty, delete
the $i^{th}$ and $k^{th}$ rows and columns of the placement on
$\lambda$ to obtain a symmetric full placement on $\hat{\lambda}$
which is counted by $\g{123}{\hat{\lambda}}{i-1}$. As all such
placements on $\hat{\lambda}$ can be produced this way, we have
$\f{123}{\lambda}{i} = \g{123}{\hat{\lambda}}{i-1}$. Similarly,
$\f{321}{\lambda}{i}$ equals the number of symmetric full
placements on $\hat{\lambda}$ which avoid $321$ and which avoid
$21$ in the $i-1$ columns $\lambda_k+1-i,\ldots,\lambda_k-1$.  By
an iterated application of Lemma~\ref{lem:inv:321slide}, this
equals $\g{321}{\hat{\lambda}}{i-1}$.  If $\hat{\lambda}$ is a
square, we may apply Lemma~\ref{lem:inv:123sq} to obtain
$\g{123}{\hat{\lambda}}{i-1} = \g{321}{\hat{\lambda}}{i-1}$ and
are thus finished. If $\hat{\lambda}$ is not a square, we show
that $\g{123}{\hat{\lambda}}{i-1} = \g{321}{\hat{\lambda}}{i-1}$
as follows.

Let $\mu$ be a non-square self-conjugate shape of length $k$; we
want to show that for $1\leq i\leq \mu_k$,
$\g{123}{\mu}{i}=\g{321}{\mu}{i}$.  If $\hat{\mu}$ is empty, then
there is a unique placement on $\mu$ and $\g{123}{\mu}{1} =
\g{321}{\mu}{1} = 1$ (and thus $\g{123}{\mu}{0} =
\g{321}{\mu}{0}$).  If $\hat{\mu}$ is non-empty, consider a
placement on $\mu$ which is counted by $\g{123}{\mu}{i}$, and let
$(j,k)$ be the position of the dot in the top row.  If $1\leq
j\leq i$, then we must have $j=1$ since columns $1,\ldots,i$ of
the placement avoid $12$. Deleting the $j^{th}$ ($1^{st}$) and
$k^{th}$ rows and columns of $\mu$, we obtain a placement on
$\hat{\mu}$ which is counted by $\g{123}{\hat{\mu}}{i-1}$ and can
obtain all such placements this way. If $j>i$, then again we
delete the $j^{th}$ and $k^{th}$ rows and columns of $\mu$, thus
obtaining those placements on $\hat{\mu}$ counted by
$\g{123}{\hat{\mu}}{j-1}$.  We thus have, for $1\leq i\leq \mu_k$,
\begin{equation}
\g{123}{\mu}{i}=\sum_{j=i}^{\mu_k} \g{123}{\hat{\mu}}{j-1}.
\label{eq:inv:g123rec}
\end{equation}
\noindent Similarly, for $1\leq i\leq \mu_k$ we also have
(invoking Lemma~\ref{lem:inv:321slide} as above)
\begin{equation}
\g{321}{\mu}{i}=\sum_{j=i}^{\mu_k} \g{321}{\hat{\mu}}{j-1}.
\label{eq:inv:g321rec}
\end{equation}
We may replace any occurrences of $\g{123}{\hat{\mu}}{0}$
($\g{321}{\hat{\mu}}{0}$) with $\g{123}{\hat{\mu}}{1}$
($\g{321}{\hat{\mu}}{1}$).  If $\mu$ is not square, we then
inductively apply this argument to $\hat{\mu}$.  Otherwise, we
apply Lemma~\ref{lem:inv:123sq} to each term in these two sums.
\end{proof}

We now complete the proof of Lemma~\ref{lem:inv:123nonsq} and
Theorem~\ref{thm:mu123} by proving Lemma~\ref{lem:inv:321slide}.
Given a $321$-avoiding placement which avoids $21$ in columns
$i,\ldots i+j-1$ but not columns $i+1,\ldots i+j$, we construct a
$321$-avoiding placement which contains $21$ in columns $i,\ldots
i+j-1$ but not in columns $i+1,\ldots i+j$.  We use a case
analysis which is more extensive than that used for general
permutations; the symmetrized version of the transformation used
in~\cite{Babson2000} works in many of the cases here but requires
adjustment in some cases.

\begin{proof}[Proof of Lemma~\ref{lem:inv:321slide}]
Note that the symmetry of the lemma imposes conditions on how the
rows of $\lambda$ are filled; the corresponding rows must also
avoid $21$ (the reflection of a $21$ pattern is also a $21$
pattern).  If $j=1$ the lemma is trivially true, so we assume that
$j\geq 2$.

Fix a symmetric full $321$-avoiding placement on $\lambda$ which
avoids $21$ in columns $i,\ldots, i+j-1$ (all of which must have
height $\lambda_1$).  Let $l$, $m$, and $r$ be the number of dots
in columns $i,\ldots i+j-1$ which are placed below, on, and above,
respectively, the diagonal of symmetry. (Since the placement
avoids $21$ in these columns, within these columns any dots which
lie above the diagonal must be to the right of all other dots,
while any dots which lie on the diagonal must be to the right of
all dots which lie below the diagonal.)\ \ Label columns $i, i+1,
\ldots, i+j-1$ as $a_1, \ldots, a_l, b_1, \ldots, b_m, c_1,
\ldots, c_r$, respectively, and denote the heights (which must
increase from left to right) of the dots in these columns by
$x_1<\cdots<x_l < y_1=b_1<\cdots<y_m=b_m < z_1<\cdots<z_r$.

We refer to the dots below, on, and above the diagonal and within
these columns as classes $\mathbf{A}$, $\mathbf{B}$, and
$\mathbf{C}$, respectively. Finally, let $w$ be the height of the
dot in column $i+j$. Figure~\ref{fig:slidesetup} shows dots in
rows and columns $i,\ldots, i+j$ of a general placement; rows and
columns $i,\ldots, i+j-1$ are marked by solid lines, while rows
and columns $i+1,\ldots, i+j$ are marked by dashed lines.  The
three classes $\mathbf{A}$, $\mathbf{B}$, and $\mathbf{C}$ of dots
are labelled and the other dots which must be present by the
symmetry of the placement are shown.  In this example, the dot in
column $i+j$ has height between the heights of two dots in class
$\mathbf{C}$.  This placement avoids $21$ in columns $i,\ldots
i+j-1$ but contains $21$ in columns $i+1, \ldots, i+j$; since it
avoids $321$ by assumption, there are various boxes (to the upper
right of the figure) which the underlying Ferrers shape must not
contain.
\begin{figure}[h]
\begin{center}
    \includegraphics[height=3.5in]{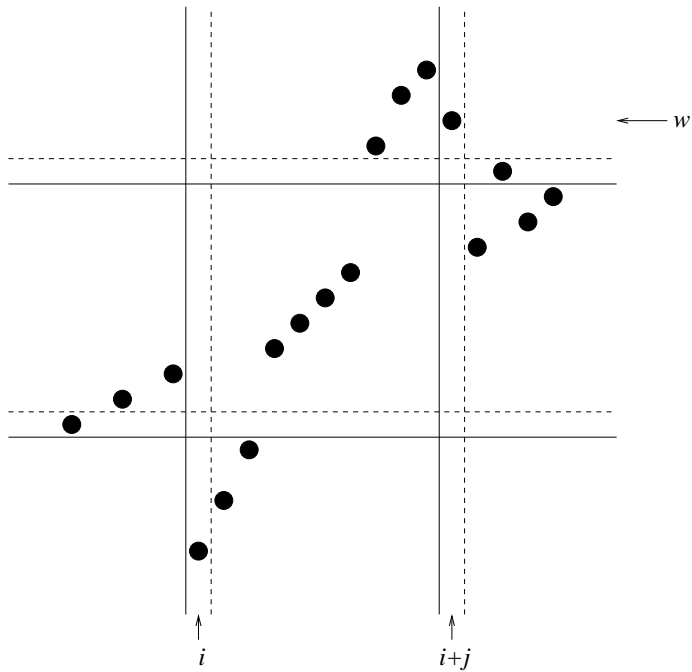}
\end{center}
\caption{The setup for the proof of
Lemma~\ref{lem:inv:321slide}.}\label{fig:slidesetup}
\end{figure}

If the placement on $\lambda$ already avoids $21$ in columns $i+1,
\ldots, i+j$, then we leave it unchanged.  If the placement
contains $21$ in these columns then we modify the placement so
that it avoids $21$ in columns $i+1, \ldots, i+j$ and contains
$21$ in columns $i, \ldots, i+j-1$.  We define our `basic
transformation' as the symmetrization of the transformation used
by Babson and West.  This suffices in most cases; we modify this
as needed below.

First of all, we will assume that $w\neq i+j$ and that if $w=c_1$
then there are no dots in class $\mathbf{B}$.  We may also assume
that exactly the $v\geq 1$ rightmost dots in columns $i, \ldots,
i+j-1$ have height greater than $w$.  Our basic transformation
consists of the following steps:
\begin{enumerate}
    \item[(i)]  Move all of the dots in columns $i,\ldots, i+j-1$ except for the
        $v^\mathrm{th}$ from the right (in column $i+j-v$, which is the
        dot in these columns with smallest height greater than $w$) $1$
        square to the right, keeping their heights unchanged.
    \item[(ii)] Move the dot in column $i+j-v$ to column $i$ (if $v=j$ this dot does not change
        position).
    \item[(iii)]    Move the dot in column $i+j$ to column $i+j-v+1$ (if
        $v=1$ this dot does not change position).
    \item[(S)]      Symmetrize these
        operations, moving the dots in the corresponding rows within their
        columns.
\end{enumerate}
Note that any dots in class $\mathbf{B}$ (on the main diagonal)
are moved $1$ square to the right by step (i) and then $1$ square
up by (S); as a result, they are still on the main diagonal in the
resulting configuration. Note also that if there is a dot in
columns $i, \ldots, i+j-1$ with height $i+j$ (in which case
$w=c_1$), this dot will be moved by both the initial operation and
the symmetrization.  In other cases (with a different position of
the dot in column $i+j$ or some dot classes being empty) the basic
transformation has similar general effects.
Figure~\ref{fig:slidebasic} shows the placement from
Figure~\ref{fig:slidesetup} using filled circles and the
transformed placement using open circles.  Arrows indicate the
general movement of dots, with filled arrowheads showing movements
due to steps (i), (ii), and (iii) and open arrowheads the
movements due to the symmetrization of these steps.  The number of
arrowheads indicates which step (or its symmetrization) is
responsible for the movement.
\begin{figure}[ht]
\begin{center}
    \includegraphics[height=3.5in]{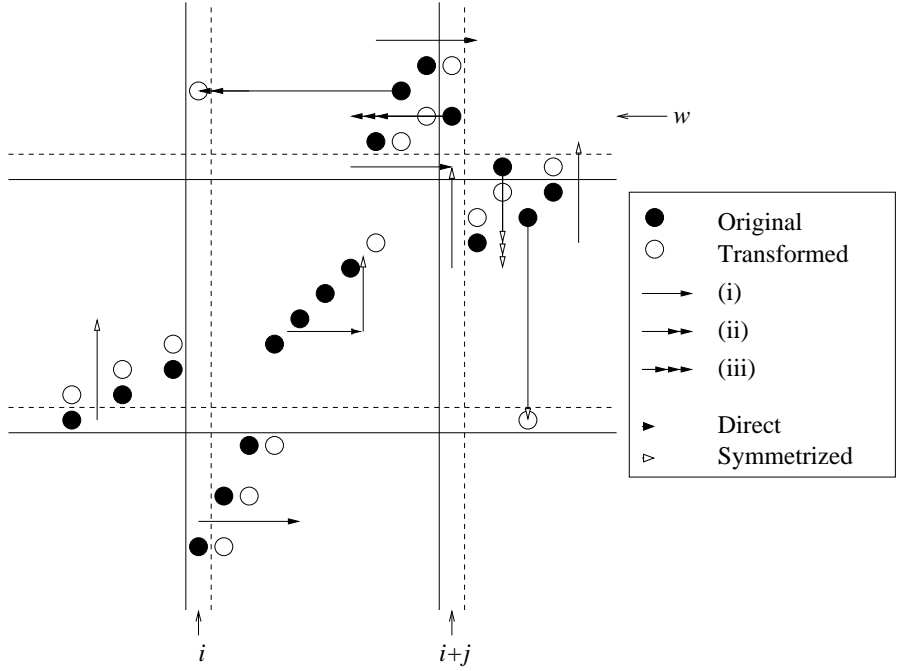}
\end{center}
\caption{The basic transformation used in proving
Lemma~\ref{lem:inv:321slide} applied to the placement in
Figure~\ref{fig:slidesetup}.}\label{fig:slidebasic}
\end{figure}

After the basic transformation has been applied, the resulting
placement avoids $21$ in columns $i+1,\ldots, i+j$ but contains it
in columns $i,\ldots, i+j-1$ since the dots which originally
formed one instance of this pattern (in columns $i+j-v$ and $i+j$)
have not changed their relative position (having moved to columns
$i$ and $i+j-v+1$) and are still contained in the rectangular
subshape of $\lambda$ formed by the leftmost $\lambda_k$ columns
(which have height $\lambda_1$).

We proceed with a case analysis on the value of $w$, modifying the
basic transformation for the case $w=c_1$ and determining the
positions of the dots in the transformed placements.  The labels
of the various subcases (\textit{e.g\/}, $\mathbf{BC}$) indicate
which of the three classes of dots are nonempty. We will then show
that this process is invertible.
\begin{description}
    \item[\textbf{I.}\ $w>\max\{x_l,y_m,z_r\}$]
        In this case, columns $i+1,\ldots,i+j$ already avoid the
        pattern $21$ and we do not modify the placement.
    \item[\textbf{II.}\ $i+j<w$]
        If this is not covered by case \textbf{I} then there is
        at least one dot in class $\mathbf{C}$ has height greater
        than $w$.  Let $d$ be the number
        of dots in class $\mathbf{C}$ with heights less than
        $w$; $w<z_1$ if $d=0$ and $z_d<w<z_{d+1}$ if $d\geq 1$.
        In each of the possible subcases --- $\mathbf{C}$,
        $\mathbf{AC}$, $\mathbf{BC}$, and $\mathbf{ABC}$ ---
        we use the basic transformation.  In columns $i+1, \ldots,
        i+j$ of the resulting configuration, $l$ dots will lie
        below the main diagonal, $m$ will lie on it, and $r$ will
        lie above it.  The dot in column $i$ will have height
        greater than exactly $d+1\geq 1$ of these dots above the
        diagonal.
    \item[\textbf{III.}\ $w=i+j$]
        In this case, the dot in column $i+j$ lies on the main
        diagonal.  If the placement is not covered by case \textbf{I}
        then class $\mathbf{C}$ must be nonempty.  In each of the
        possible subcases --- $\mathbf{C}$, $\mathbf{AC}$,
        $\mathbf{BC}$, and $\mathbf{ABC}$ --- we may use the basic
        transformation.  Note that the dot initially in
        $(i+j,i+j)$ is moved to $(c_1+1,i+j)$ by step (iii) of the
        basic transformation and then to $(c_1+1,c_1+1)$ by the
        symmetrization step (S).

        In columns $i+1, \ldots, i+j$ of the resulting
        configuration $r-1\geq 0$ dots will lie above the main
        diagonal, $m+1\geq 1$ will lie on it, and $l\geq 0$ will
        lie below it.  The height of the dot in column $i$ will be
        greater than the heights of all of these dots below and on
        the diagonal but less than the heights of all of these
        dots above the diagonal.
    \item[\textbf{IV.}\ $w=c_1, m=0$]
        In this case, the point in column $i+j$ is the
        reflection across the main diagonal of the leftmost dot in
        class $\mathbf{C}$ (which is thus nonempty).  Since there
        are no dots in class $\mathbf{B}$, we may use the basic
        transformation for the two possible subcases
        ($\mathbf{C}$ and $\mathbf{AC}$).

        In columns $i+1, \ldots, i+j$ of the resulting
        configuration, $r-1\geq 0$ dots lie above the main
        diagonal, $0$ lie on it, and $l+1 \geq 1$ lie below it.
        The dot in column $i$ will have height greater than the
        heights of these dots which lie below the diagonal and
        less than the heights of these dots which lie above the
        diagonal.  Note that the dots originally in $(i+j,c_1)$
        and $(c_1,i+j)$ are moved to $(c_1+1,i)$ and $(i,c_1+1)$
        by the combined effects of steps (i) and (S) of the basic
        transformation.
    \item[\textbf{V.}\ $w=c_1, m>0$]
        In this case, the point in column $i+j$ is the
        reflection across the main diagonal of the leftmost dot in
        class $\mathbf{C}$ (which is thus nonempty).  As
        class $\mathbf{B}$ is nonempty we need to modify the basic
        transformation as follows.  Move the dots in columns $c_2,
        \ldots, c_r$ and any dots in columns $a_1, \ldots, a_l$
        one box to the right and symmetrize this.  (This is step (i) and its
        symmetrization from the basic transformation.)\ \ Combine
        the $2$ dots at $(i+j,c_1)$ and $(c_1,i+j)$ into
        a single dot at $(c_1+1,c_1+1)$.  Split the dot in
        $(b_1,b_1)$ into $2$ dots at $(i,b_1+1)$ and $(b_1+1,i)$.
        Move any dots in columns $b_2, \ldots, b_m$ one box up and
        right, keeping them on the main diagonal.

        In columns $i+1, \ldots, i+j$ of the resulting
        configuration, $r-1\geq 0$ dots lie above the main
        diagonal, $m\geq 1$ lie on it, and $l+1\geq 1$ lie below
        it.  The dot in column $i$ has height greater than exactly
        $l+1$ of these dots.
    \item[\textbf{VI.}\ $w<i$]    In this case, class $\mathbf{B}$
    must be empty since
    the dots $(w,i+j)$, $(y_1,b_1)$, and $(i+j,w)$ would form a
    $321$ pattern.  Any dots in class $\mathbf{A}$ must have height less
    than $w$, since if $x_l>w$ then the dots $(w,i+j)$, $(x_l,a_l)$,
    and $(i+j,w)$ would form a $321$ pattern.  We use the basic
    transformation for the subcases --- $\mathbf{C}$ and
    $\mathbf{AC}$ --- not covered by case \textbf{I} above.

    In columns $i+j, \ldots, i+j$ of the resulting configuration,
    $r-1\geq 0$ dots lie above the main diagonal, $0$ lie on it,
    and $l+1\geq 1$ lie below it.  The dot in column $i$ will have
    height greater than $i+j$ and less than the heights of any of
    these dots above the diagonal.
\end{description}

The basic transformation is invertible; the modified version used
in case \textbf{V} is also invertible and has an image which does
not meet that of the basic transformation.

We now see that the possible arrangements of $321$-avoiding
placements which avoid $21$ in columns $i+1, \ldots, i+j$ are all
in the image of one of the transformations above.  Let $u$ be the
height of the dot in column $i$.  If $u$ is less than the heights
of all the dots in columns $i+1, \ldots, i+j$ we do nothing (case
\textbf{I} above); this covers all $u\leq i$.  If $i<u\leq i+j$
then the arrangement is in the image of either case \textbf{IV}
(if there are no dots on the diagonal in these columns) or case
\textbf{V} (if there are).  If $i+j<u$ and there are no dots in
these columns above the diagonal with height less than $u$ then
the arrangement is in the image of either case \textbf{VI} (if
there are no dots on the diagonal in these columns) or case
\textbf{III} (if there are).  Finally, if $i+j<u$ and there is at
least one dot in these columns with height less than $u$, then the
arrangement is in the image of case \textbf{II} above.
\end{proof}

\section{Classifying avoided patterns}
\label{sec:conc}

Corollary~\ref{cor:ex12} implies many of the previously known
$\ieq$-equivalences between symmetry classes in $\S{4}$. The most
notable of these is $1234\ieq 2143$, originally conjectured by
Guibert~\cite{GuibertThesis} and proved more recently (using
different arguments than we use) by Guibert, Pergola, and
Pinzani~\cite{Guibert2001}.
\begin{cor}[Guibert, Pergola, Pinzani \cite{Guibert2001}]
\begin{equation*}
1234\ieq 2134\ieq 2143
\end{equation*}
\end{cor}
\begin{proof}
Apply Corollary~\ref{cor:ex12} to both members of the symmetry
class $\{1243,2134\}$.
\end{proof}

Corollary~\ref{cor:ex123} implies an affirmative answer to a
conjecture of Guibert (Conjecture~5.3 of~\cite{Guibert2001},
originally from~\cite{GuibertThesis}) as follows.
\begin{cor}
\begin{equation*}
1234\ieq 3214
\end{equation*}
\label{cor:inv:3214}
\end{cor}
\noindent Corollary~\ref{cor:inv:3214} completes the
classification of $\S{4}$ according to $\ieq$-e\-quiv\-a\-lence;
we show this classification in Table~\ref{tab:inv:new4patt}.  The
leftmost column contains the symmetry classes of $\S{4}$, denoted
by braces $\{\cdot\}$. Double lines mark the partition of the set
of symmetry classes into cardinality classes. The remaining
columns give $\iav{n}{\tau}$ for $\tau$ in each of the cardinality
classes and $5\leq n\leq 11$.  As these sequences differ from each
other, there can be no other $\ieq$-equivalences in $\S{4}$.
\begin{table}[h]
\begin{center}
\begin{tabular}{|l|r|r|r|r|r|r|r|}\hline
$\tau$ & $\iav{5}{\tau}$ & $\iav{6}{\tau}$ & $\iav{7}{\tau}$ &
$\iav{8}{\tau}$ & $\iav{9}{\tau}$ & $\iav{10}{\tau}$ &
$\iav{11}{\tau}$ \\ \hline\hline
$\{1324\}$              & 21 & 51 & 126 & 321 & 820 & 2160 & 5654 \\
\hline\hline
$\{1234\}$              & 21 & 51 & 127 & 323 & 835 & 2188 & 5798\\ 
$\{1243, 2134\}$        & & & & & & &\\ 
$\{1432, 3214\}$        & & & & & & &\\
$\{2143\}$              & & & & & & &\\
$\{3412\}$              & & & & & & &\\
$\{4321\}$              & & & & & & &\\ \hline\hline
$\{4231\}$              & 21 & 51 & 128 & 327 & 858 & 2272 & 6146\\
\hline\hline
$\{2431, 4132,$         & 24 & 62 & 154 & 396 & 992 & 2536 & 6376\\ 
$\quad 3241, 4213\}$ & & & & & & &\\ 
\hline\hline
$\{1342, 1423,$         & 24 & 62 & 156 & 406 & 1040 & 2714 & 7012\\
$\quad 2314, 3124\}$    & & & & & & &\\
\hline\hline
$\{2341, 4123\}$        & 25 & 66 & 170 & 441 & 1124 & 2870 & 7273\\
\hline\hline
$\{3421, 4312\}$        & 25 & 66 & 173 & 460 & 1218 & 3240 & 8602\\
\hline\hline
$\{2413, 3142\}$        & 24 & 64 & 166 & 456 & 1234 & 3454 & 9600\\ 
\hline
\end{tabular}
\end{center}
\caption{The completed classification of $\S{4}$ according to
pattern avoidance by involutions.  Values of $\iav{n}{\tau}$ for
$\tau\in\S{4}$ and $5\leq n\leq 11$.} \label{tab:inv:new4patt}
\end{table}

In light of the results of Sections~\ref{sec:12tau}
and~\ref{sec:123tau}, it is natural to conjecture that similar
theorems hold for $12\ldots k$ and $k\ldots 21$.
\begin{conj}
For every $k$, the prefixes $12\ldots k$ and $k\ldots 21$ may be
exchanged. \label{conj:exidk}
\end{conj}
\noindent  We note that $1234\ieq 3412$ is the only
$\ieq$-equivalence in $\S{4}$ that does not follow from
Conjecture~\ref{conj:exidk}.  More generally, we also make the
following conjecture.
\begin{conj}
Let $\mu_{sym}(\sigma)$ be the number of symmetric full placements
on the shape $\mu$ which avoid the pattern $\sigma$.  For every
$k$ and self-conjugate shape $\mu$,
\begin{equation*}
\mu_{sym}(12\ldots k) = \mu_{sym}(k\ldots 21)
\end{equation*}
\label{conj:muidk}
\end{conj}
\noindent As noted above, Backelin, West, and
Xin~\cite{Backelin2003} have recently proved the analogue of this
conjecture for pattern avoidance by general permutations.

Corollary~\ref{cor:ex12} and~\ref{cor:ex123} also imply
(apparently new) $\ieq$-equivalences for patterns in $\S{5}$.
Numerical results indicate that among the non-involutions in
$\S{5}$ there is only one possible $\ieq$-equivalence, while
Corollary~\ref{cor:ex12} implies that this does indeed hold.
\begin{cor}
\begin{equation*}
12453\ieq 21453
\end{equation*}
\end{cor}
Among the involutions of length $5$, one cardinality class
contains at most the symmetry classes of $12435$ and $21435$; this
$\ieq$-equivalence also holds.
\begin{cor}
\begin{equation*}
12435\ieq 21435
\end{equation*}
\end{cor}
\noindent Numerical results suggest that many symmetry classes may
belong to the same cardinality class as $\{12345\}$.
Corollaries~\ref{cor:ex12} and~\ref{cor:ex123} imply some of these
possible $\ieq$-equivalences.
\begin{cor}
\begin{equation*}
21543\ieq 12543\ieq 12345\ieq 12354\ieq 21354
\label{cor:s5idclass}
\end{equation*}
\end{cor}
\noindent Recall that even without Conjecture~\ref{conj:exidk}, we
have $12\ldots k\ieq k\ldots 21$ by the symmetry of the RSK
algorithm.

Table~\ref{tab:inv:new5pattinvols} displays the involutions in
$\S{5}$ using the same format as Table~\ref{tab:inv:new4patt},
with the addition of single lines to separate symmetry classes
with identical numerical results but which have not yet been
proven to be in the same cardinality class.
Table~\ref{tab:inv:new5pattnoninvols} shows the non-involutions in
$\S{5}$ in this format.
\begin{table}
\begin{center}
\begin{tabular}{|l|r|r|r|r|r|r|}\hline
$\tau$ & $\iav{6}{\tau}$ & $\iav{7}{\tau}$ & $\iav{8}{\tau}$ &
$\iav{9}{\tau}$ & $\iav{10}{\tau}$ & $\iav{11}{\tau}$\\
\hline\hline
$\{35142, 42513\}$          & 70 & 195 & 582 & 1725 & 5355 & 16510\\
\hline\hline
$\{14325\}$                 & 70 & 196 & 587 & 1757 & 5504 & 17220\\
\hline\hline
$\{12435, 13245\}$          & 70 & 196 & 587 & 1759 & 5512 & 17290\\
$\{13254, 21435\}$          & & & & & & \\ \hline\hline
$\{12345\}$                 & 70 & 196 & 588 & 1764 & 5544 & 17424\\
$\{54321\}$                 & & & & & & \\
$\{12354, 21345\}$          & & & & & & \\
$\{12543, 32145\}$          & & & & & & \\
$\{21354\}$                 & & & & & & \\
$\{21543, 32154\}$          & & & & & & \\
\hline
$\{15432, 43215\}$          & 70 & 196 & 588 & 1764 & 5544 & 17424\\
\hline
$\{45312\}$                 & 70 & 196 & 588 & 1764 & 5544 & 17424\\
\hline\hline
$\{52431, 53241\}$          & 70 & 196 & 588 & 1764 & 5544 & 17426\\
\hline\hline
$\{52341\}$                 & 70 & 196 & 589 & 1773 & 5604 & 17768\\
\hline\hline
$\{14523, 34125\}$          & 70 & 197 & 592 & 1791 & 5644 & 17900\\
\hline\hline
$\{15342, 42315\}$          & 70 & 197 & 593 & 1797 & 5685 & 18101\\
\hline
\end{tabular}
\end{center}
\caption{Values of $\iav{n}{\tau}$ for involutions $\tau\in\S{5}$
and $6\leq n\leq 11$.} \label{tab:inv:new5pattinvols}
\end{table}
\begin{table}
\begin{center}
\begin{tabular}{|l|r|r|r|r|r|r|}\hline
$\tau$ & $\iav{6}{\tau}$ & $\iav{7}{\tau}$ & $\iav{8}{\tau}$ &
$\iav{9}{\tau}$ & $\iav{10}{\tau}$ & $\iav{11}{\tau}$\\
\hline\hline $\{13542, 42135, 15243, 32415\}$    & 74 & 214 & 644
& 1945 & 6004 & 18526\\ \hline\hline $\{13425, 14235\}$
& 74 & 214 & 647 & 1959 & 6107 & 18952\\ \hline\hline $\{14352,
41325, 15324, 24315\}$    & 74 & 215 & 649 & 1975 & 6126 & 19057\\
\hline\hline $\{25431, 53214, 43251, 51432\}$    & 74 & 216 & 654
& 2002 & 6223 & 19425\\ \hline\hline $\{45231, 53412\}$
& 74 & 216 & 656 & 2020 & 6342 & 20072\\ \hline\hline
$\{12453, 31245, 12534, 23145\}$    & 74 & 216 & 656 & 2022 & 6362 & 20212\\
$\{21453, 31254, 21534, 23154\}$    & & & & & & \\ \hline\hline
$\{25143, 32514, 31542, 42153\}$    & 74 & 216 & 658 & 2033 & 6434
& 20538\\ \hline\hline $\{13452, 41235, 15234, 23415\}$    & 75 &
220 & 674 & 2067 & 6463 & 20150\\ \hline\hline $\{13524, 24135,
14253, 31425\}$    & 74 & 217 & 664 & 2068 & 6598 & 21269\\
\hline\hline $\{14532, 43125, 15423, 34215\}$    & 75 & 220 & 677
& 2090 & 6609 & 20880\\ \hline\hline $\{32541, 52143\}$
& 75 & 221 & 679 & 2096 & 6577 & 20630\\ \hline\hline $\{25341,
52314, 42351, 51342\}$    & 75 & 221 & 680 & 2103 & 6617 & 20808\\
\hline\hline $\{35241, 52413, 42531, 53142\}$    & 75 & 220 & 680
& 2111 & 6745 & 21567\\ \hline\hline $\{24513, 35124, 34152,
41523\}$    & 75 & 221 & 682 & 2122 & 6752& 21569\\ \hline\hline
$\{53421, 54231\}$                  & 74 & 218 & 672 & 2126 & 6908
& 22877\\ \hline\hline $\{24351, 51324\}$                  & 75 &
222 & 687 & 2136 & 6735 & 21093\\ \hline\hline $\{23541, 52134,
32451, 51243\}$    & 75 & 222 & 687 & 2137 & 6737 & 21132\\
\hline\hline $\{45321, 54312\}$                  & 75 & 222 & 688
& 2156 & 6892 & 22128 \\ \hline\hline $\{23514, 25134, 31452,
41253\}$    & 75 & 222 & 688 & 2159 & 6923 & 22358\\ \hline\hline
$\{35412, 45213, 43512, 45132\}$    & 75 & 222 & 689 & 2168 & 6981
& 22676\\ \hline\hline $\{25413, 35214, 41532, 43152\}$    & 75 &
222 & 690 & 2172 & 7004 & 22731\\ \hline\hline $\{23451, 51234\}$
& 75 & 223 & 694 & 2183 & 6958 & 22127\\ \hline\hline $\{35421,
54213, 43521, 54132\}$    & 75 & 223 & 696 & 2209 & 7177 & 23533\\
\hline\hline $\{24153, 31524\}$                  & 75 & 224 & 701
& 2240 & 7315 & 24190\\ \hline\hline $\{24531, 53124, 34251,
51423\}$    & 76 & 227 & 715 & 2257 & 7269 & 23254\\ \hline\hline
$\{34512, 45123\}$                  & 75 & 224 & 705 & 2273 & 7538
& 25418\\ \hline\hline $\{25314, 41352\}$                  & 76 &
228 & 722 & 2302 & 7514 & 24530\\ \hline\hline
$\{34521, 54123\}$                  & 76 & 230 & 732 & 2364 & 7764 & 25596\\ 
\hline
\end{tabular}
\end{center}
\caption{Values of $I_n(\tau)$ for non-involutions $\tau\in\S{5}$
and $6\leq n\leq 11$.} \label{tab:inv:new5pattnoninvols}
\end{table}

We conjecture that both of the possible equivalences consistent
with Table~\ref{tab:inv:new5pattinvols} do hold.
\begin{conj}
\begin{equation*}
12345\ieq 43215
\end{equation*}
\end{conj}
\begin{conj}
\begin{equation*}
12345\ieq 45312
\end{equation*}
\end{conj}
\noindent The first of these follows from
Conjecture~\ref{conj:exidk}, while the second is the only
$\ieq$-equivalence which is consistent with the numerical results
presented here but which would not follow from
Conjecture~\ref{conj:exidk}.

Table~\ref{tab:idclass6inv} uses a somewhat different format to
present just those involutions $\tau\in\S{6}$ for which
$\iav{n}{\tau}=\iav{n}{123456}$ for $7\leq n\leq 11$.  Single
lines now separate symmetry classes whose $\ieq$-equivalence
follows from Conjecture~\ref{conj:exidk}, but which are not yet
known to be in the same cardinality class, while double lines
separate the other symmetry classes whose numerical results agree
with those of $123456$.  Symmetry classes which are not separated
by any lines are known to be in the same cardinality class.
\begin{table}
\begin{center}
\begin{tabular}{|l|r|r|r|r|r|}\hline
$\tau$ &  $\iav{7}{\tau}$ & $\iav{8}{\tau}$ &
$\iav{9}{\tau}$ & $\iav{10}{\tau}$ & $\iav{11}{\tau}$\\
\hline\hline
$\{123456\}$                                    & 225 & 715 & 2347 & 7990 & 27908\\
$\{123465, 213456\}$                            &  &  &  &  & \\
$\{123654, 321456\}$                            &  &  &  &  & \\
$\{213465\}$                                    &  &  &  &  & \\
$\{213654, 321465\}$                            &  &  &  &  & \\
$\{321654\}$                                    &  &  &  &  & \\
$\{654321\}$                                    &  &  &  &  & \\
\hline
$\{126543, 432156\}$                            & 225 & 715 & 2347 & 7990 & 27908\\
$\{216543, 432165\}$                            &  &  &  &  & \\
\hline $\{165432, 543216\}$                            & 225 & 715
& 2347 & 7990 & 27908\\ \hline\hline $\{456123\}$
& 225 & 715 & 2347 & 7990 & 27908\\ \hline\hline
$\{564312\}$                                    & 225 & 715 & 2347 & 7990 & 27908\\ \hline
\end{tabular}
\end{center}
\caption{Values of $\iav{n}{\tau}$ for those $\tau\in\S{6}$ which
may be $\ieq$-equivalent to $123456\in\S{6}$ and $7\leq n\leq
11$.\label{tab:idclass6inv}}
\end{table}

A natural question is whether there are any $\ieq$-equivalences
between $123456$ and other patterns which do not follow from
Conjecture~\ref{conj:exidk}.
\begin{q}
Does either $123456\ieq 456123$ or $123456\ieq 564312$ hold?
\end{q}
\noindent These resemble the known $\ieq$-e\-quiv\-a\-lence
$1234\ieq 3412$ and the conjectured $\ieq$-e\-quiv\-a\-lence
$12345\ieq 45312$, the only possible $\ieq$-equivalences in
$\S{4}$ and $\S{5}$ which do not follow from
Conjecture~\ref{conj:exidk}. It is natural to ask whether these
follow from some other general theorem.  Such a general result
could not be stated in terms of prefix exchanging as
$12345\not\ieq 34125$ and $123456\not\ieq 453126$ (although not
shown in Table~\ref{tab:idclass6inv}, $\iav{8}{453126}=716> 715 =
\iav{8}{123456}$).
\begin{q}
Is there a general theorem which implies $1234\ieq 3412$,
$12345\ieq 45312$, and one or both of $123456\ieq 456123$ and
$123456\ieq 564312$?
\end{q}

Stankova and West~\cite{Stankova2002} proved that $231$ and $312$
are shape-Wilf-equivalent. However, we cannot have
$\mu_{sym}(231)=\mu_{sym}(312)$ for every self-conjugate shape
$\mu$ since $\iav{10}{231564} = 8990 < 8991 = \iav{10}{312564}$.
Stankova and West also noted a `sporadic' case of Wilf-equivalence
(between the $\seq$-symmetry classes of $1342$ and $3142$) which
does not follow from any known instance of shape-Wilf-equivalence.
It is interesting to see that this Wilf-equivalence breaks when we
pass to $\ieq$-cardinality classes; neither of the two
$\ieq$-symmetry classes contained in the $\seq$-symmetry class of
$1342$ are $\ieq$-equivalent to the $\ieq$-symmetry class
containing $3142$ (which in also the $\seq$-symmetry class of
$3142$).

There are also asymptotic questions which can be asked for pattern
avoidance by involutions.  It seems reasonable to ask, as Stankova
and West have conjectured for the general permutation case,
whether the cardinality classes can be ordered asymptotically. For
relatively small values of $n$, it is reasonable to expect that if
$\sigma$ is an involution and $\pi$ is not then
$\iav{n}{\sigma}<\iav{n}{\pi}$.  However, comparing the growth
rates of, \textit{e.g.\/}, $\iav{n}{15342}$ and $\iav{n}{13542}$,
suggests that it may be too optimistic to expect this to hold for
all $n$. Finally, it would be of interest to find the asymptotic
growth of $\iav{n}{\tau}$ for those $\tau\not\ieq 12\ldots k$.

\section*{Acknowledgements}
We are grateful to Herb Wilf for helpful discussions about this
work, the organizers of and participants in the Permutation
Patterns 2003 conference for stimulating interactions about
current pattern work, and Andre Scedrov for arranging research
support.  We also thank the anonymous referee for a number of
useful suggestions.

\providecommand{\bysame}{\leavevmode\hbox
to3em{\hrulefill}\thinspace}
\providecommand{\MR}{\relax\ifhmode\unskip\space\fi MR }
\providecommand{\MRhref}[2]{%
  \href{http://www.ams.org/mathscinet-getitem?mr=#1}{#2}
} \providecommand{\href}[2]{#2}

\end{document}